\documentclass[11pt]{article}

\usepackage[T1]{fontenc}
\usepackage[utf8]{inputenc}
\usepackage{lmodern}
\usepackage[letterpaper,margin=1in]{geometry}
\usepackage{amsmath,amssymb,amsthm,mathtools}
\usepackage{microtype}
\usepackage{booktabs}
\usepackage{xurl}
\usepackage{needspace}
\usepackage[dvipsnames]{xcolor}
\usepackage[colorlinks=true,linkcolor=MidnightBlue,
            citecolor=MidnightBlue,urlcolor=MidnightBlue]{hyperref}
\hypersetup{
  pdftitle={Packed Passengers in Binary Seating Arrangements},
  pdfauthor={Perplexity Computer},
  pdfsubject={Enumerative combinatorics of overlapping triples and grid neighborhoods},
  pdfkeywords={binary strings, overlapping patterns, generating functions, Fibonacci numbers}
}

\newtheorem{theorem}{Theorem}[section]
\newtheorem{proposition}[theorem]{Proposition}
\newtheorem{corollary}[theorem]{Corollary}
\theoremstyle{definition}
\newtheorem{definition}[theorem]{Definition}
\theoremstyle{remark}

\newtheorem{example}[theorem]{Example}
\newcommand{\Fcal}{\mathcal{F}}
\newcommand{\Hcal}{\mathcal{H}}
\newcommand{\Prob}{\mathbb{P}}
\newcommand{\E}{\mathbb{E}}
\numberwithin{equation}{section}
\allowdisplaybreaks[1]

\title{On Packed Passenger Permutations}
\author{Erik Dahlen}
\date{}

\begin{document}
\maketitle

\begin{abstract}
An occupied seat in a binary row is called packed if both of its immediate
neighbors are occupied. We give a self-contained combinatorial account of
this statistic, equivalently the number of overlapping occurrences of
\(111\). An occupied-block decomposition yields a coefficient formula at
fixed occupancy and a direct proof of the sharp minimum
\(\max\{0,3m-2n-2\}\) for \(m\) passengers in \(n\) seats, including
arrangements attaining the bound. The same decomposition gives a rational
generating function, recurrences, expectations, and avoidance asymptotics.
Full rows, in which every empty seat would become packed upon being filled,
are related to independent sets of a path. A neighborhood-union expansion
also gives exact counts on rectangular grids. The emphasis is a unified
elementary presentation and explicit extremal constructions, rather than
a new general enumeration method.
\end{abstract}

\section{Introduction and notation}
\label{sec:intro}

Consider a row of seats, each either occupied or empty. A passenger is
packed when the seats immediately to the left and right are both occupied;
passengers at the ends of the row are never packed. Although the seating
interpretation is convenient, the underlying objects are binary strings,
not permutations of distinct passengers.

The packed-passenger statistic counts overlapping occurrences of the word
\(111\). General methods for string overlaps and pattern avoidance were
developed by Guibas and Odlyzko \cite{GO}.\footnote{Journal article:
\url{https://doi.org/10.1016/0097-3165(81)90005-4}.}
The weighted enumeration of word occurrences, including the simultaneous
recording of letter counts, is treated by the Goulden--Jackson cluster
method; see Noonan and Zeilberger \cite{NZ}.\footnote{Full text:
\url{https://arxiv.org/html/math/9806036}.}
More directly, Inoue and Aki give a first-failure recurrence and
success-block generating functions for overlapping success runs
\cite[Section~2]{IA05}.\footnote{Full text:
\url{https://www.ism.ac.jp/editsec/aism/pdf/057_2_0353.pdf}.}
Their subsequent treatment includes run distributions conditional on
fixed success counts \cite[Section~4]{IA07}.\footnote{Full text:
\url{https://www.ism.ac.jp/editsec/aism/pdf/059_3_0577.pdf}.}
The row distributions studied here are the length-three Type~III
(fully overlapping) specializations of that theory.

Our contribution is a unified elementary exposition of this seating
model, centered on the occupied-block decomposition and the direct
extremal calculation
\[
 P(m,n)=\max\{0,3m-2n-2\}.
\]
We give explicit arrangements attaining this bound, relate full rows
to independent sets of a path, and work out a rectangular-grid example
through neighborhood unions. We also give unique combinatorial solutions, using coefficients of generating functions to find this minimum and count arrangements with $k$ amount of 3-runs.

Section~\ref{sec:fixed} treats a fixed number of passengers, and
Section~\ref{sec:minimum} determines the minimum packed count.
Section~\ref{sec:unrestricted} sums over all passenger counts.
Section~\ref{sec:probability} develops probabilistic consequences,
Section~\ref{sec:full} counts full rows, and
Section~\ref{sec:grid} considers rectangular grids.

For \(n\geq0\), write \(s=s_1\cdots s_n\in\{0,1\}^n\), where \(s_i=1\)
means that seat \(i\) is occupied. Set
\[
 m(s)=\sum_{i=1}^n s_i,\qquad
 S_{m,n}=\{s\in\{0,1\}^n:m(s)=m\},\qquad
 |S_{m,n}|=\binom{n}{m}.
\]
The empty string is included when \(n=0\).

\Needspace{5\baselineskip}
\begin{definition}[Packed count]
The number of packed passengers in \(s\) is
\begin{equation}\label{eq:p}
 p(s)=\sum_{i=2}^{n-1}s_{i-1}s_is_{i+1}.
\end{equation}
An empty sum is zero, so \(p(s)=0\) for \(n<3\).
Occurrences may overlap: for example, \(p(1111)=2\).
\end{definition}

For integers \(0\leq m\leq n\), define
\begin{align}
 G_{m,n}(z)&=\sum_{s\in S_{m,n}}z^{p(s)},\label{eq:Gdef}\\
 P(m,n)&=\min_{s\in S_{m,n}}p(s).\label{eq:Pdef}
\end{align}
Thus \([z^k]G_{m,n}(z)\) counts strings with \(m\) occupied seats and \(k\)
packed passengers. We use \([z^k]\) for coefficient extraction.
The unrestricted polynomial and its ordinary generating function are
\begin{equation}\label{eq:Hdef}
 H_n(z)=\sum_{m=0}^n G_{m,n}(z),\qquad
 \Hcal(x,z)=\sum_{n\geq0}H_n(z)x^n.
\end{equation}
Generating-function identities are interpreted as formal power-series
identities, unless a numerical specialization is explicitly stated.

\section{The distribution of a fixed passenger count}
\label{sec:fixed}

Let \(f(k)=\max\{0,k-2\}\). An occupied block of length \(k\) contributes
exactly \(f(k)\) packed passengers.

\subsection{The first-block recurrence}

\begin{proposition}\label{prop:Grec}
For \(0\leq m<n\),
\begin{equation}\label{eq:Grec}
 G_{m,n}(z)=\sum_{k=0}^m z^{f(k)}G_{m-k,n-k-1}(z).
\end{equation}
The initial cases are \(G_{r,r}(z)=z^{f(r)}\) for \(r\geq0\).
\end{proposition}

\begin{proof}
There is at least one zero. Cutting at the first zero gives the unique
decomposition \(s=1^k0s'\), with \(0\leq k\leq m\), where
\(s'\in S_{m-k,n-k-1}\). No occupied triple crosses the separating zero,
so \(p(s)=f(k)+p(s')\). Summing the resulting weights proves the
recurrence. If \(m=n=r\), the only string is \(1^r\).
\end{proof}

\begin{example}\label{ex:G56}
With five passengers in six seats, there is one zero, and
\[
 G_{5,6}(z)=\sum_{k=0}^5 z^{f(k)+f(5-k)}
          =2z+2z^2+2z^3.
\]
In particular, \(G_{5,6}(1)=6=\binom65\).
\end{example}

\subsection{Occupied-block generating functions}

The generating function for a possibly empty occupied block is
\begin{equation}\label{eq:B}
 B(t,z)=\sum_{k\geq0}t^kz^{f(k)}
       =1+t+t^2+\frac{t^3z}{1-tz}.
\end{equation}
Here \(t\) records occupied seats and \(z\) records packed passengers.
Empty blocks allow consecutive zeros and zeros at either endpoint.

\begin{theorem}[Fixed-count distribution]\label{thm:Ggf}
For integers \(0\leq m\leq n\),
\begin{equation}\label{eq:Ggf}
 G_{m,n}(z)=[t^m]\left(1+t+t^2+\frac{t^3z}{1-tz}\right)^{n-m+1}.
\end{equation}
\end{theorem}

\begin{proof}
Put \(q=n-m\). Every string has the unique form
\[
 1^{k_0}0\,1^{k_1}0\cdots0\,1^{k_q},
 \qquad k_i\geq0,\qquad k_0+\cdots+k_q=m.
\]
The contribution to the packed count is \(\sum_{i=0}^q f(k_i)\).
Consequently,
\begin{equation}\label{eq:Gblocks}
 G_{m,n}(z)=
 \sum_{\substack{k_0+\cdots+k_q=m\\k_i\geq0}}
 z^{\sum_{i=0}^q f(k_i)}
 =[t^m]B(t,z)^{q+1}.
\end{equation}
\end{proof}

Equivalently, fixing the number \(q\) of zeros gives
\begin{equation}\label{eq:Fq}
 \Fcal_q(t,z):=\sum_{m\geq0}G_{m,m+q}(z)t^m=B(t,z)^{q+1}.
\end{equation}
This also follows from Proposition~\ref{prop:Grec}, which gives
\(\Fcal_q=B\Fcal_{q-1}\) for \(q\geq1\), with \(\Fcal_0=B\).
At \(z=1\), \(B(t,1)=(1-t)^{-1}\), so
\eqref{eq:Ggf} recovers \(G_{m,n}(1)=\binom nm\).

For comparison with the fixed-composition literature, record the
occupied-seat count by \(u\) and the total length by \(x\). Then
\[
 \sum_{n\geq0}\sum_{m=0}^n G_{m,n}(z)u^mx^n
 =\frac{B(ux,z)}{1-xB(ux,z)}.
\]
This is the one-success-type, length-three Type~III specialization
of Inoue and Aki's conditional generating function
\cite[Theorem~4.2]{IA07}.\footnote{Theorem~4.2, equation~(18):
\url{https://www.ism.ac.jp/editsec/aism/pdf/059_3_0577.pdf}.}
The coefficient formula above follows by fixing the number of
separating zeros.

\section{The minimum number of packed passengers}
\label{sec:minimum}

All coefficients of \(G_{m,n}\) are nonnegative, and therefore
\begin{equation}\label{eq:Pminexp}
 P(m,n)=\min\{k\geq0:[z^k]G_{m,n}(z)>0\}.
\end{equation}
The block decomposition makes this minimum explicit.

\begin{theorem}[Minimum packed count]\label{thm:P}
For integers \(0\leq m\leq n\),
\begin{equation}\label{eq:Pformula}
 P(m,n)=\max\{0,3m-2n-2\}.
\end{equation}
\end{theorem}

\begin{proof}
Let \(q=n-m\) and use the block lengths in \eqref{eq:Gblocks}. Since
\[
 \sum_{i=0}^q f(k_i)
 =m-\sum_{i=0}^q\min\{k_i,2\}
 \geq m-2(q+1),
\]
and the left side is nonnegative, this gives the required lower bound.
If \(m\leq2(q+1)\), distribute the \(m\) occupied seats among the \(q+1\)
blocks with each block of length at most \(2\). This attains zero.
If \(m>2(q+1)\), choose
\[
 k_0=m-2q,\qquad k_1=\cdots=k_q=2.
\]
Then \(k_0>2\), and the packed count is \(m-2(q+1)\).
For \(q=0\), this construction consists only of \(k_0=m\).
In both cases the bound is attained. Substituting \(q=n-m\) completes
the proof.
\end{proof}

\subsection{A local identity}

An alternative lower-bound argument uses adjacent pairs and isolated
occupied seats. Set \(s_0=s_{n+1}=0\), and write
\[
 A(s)=\sum_{i=1}^{n-1}s_is_{i+1},\qquad
 U(s)=\sum_{i=1}^n s_i(1-s_{i-1})(1-s_{i+1}).
\]
The zero padding is used only in these expressions; it does not change
the seating model.

\begin{proposition}\label{prop:local}
For every binary string \(s\),
\begin{equation}\label{eq:localidentity}
 p(s)=2A(s)-m(s)+U(s).
\end{equation}
\end{proposition}

\begin{proof}
Expanding the expression for \(U(s)\), the two adjacent-pair sums each
equal \(A(s)\), and the triple sum equals \(p(s)\). Hence
\(U(s)=m(s)-2A(s)+p(s)\).
\end{proof}

For a string with \(m\) occupied seats and \(q=n-m\) empty seats,
let \(b\) be the number of nonempty occupied blocks. Then \(A(s)=m-b\).
Since \(b\leq q+1\) and \(U(s)\geq0\), \eqref{eq:localidentity} gives
\(p(s)\geq3m-2n-2\). Together with \(p(s)\geq0\), this is the lower
bound in Theorem~\ref{thm:P}; its proof supplies the attaining strings.

\section{All passenger counts}
\label{sec:unrestricted}

By definition, \([z^k]H_n(z)\) counts all binary strings of length \(n\)
with exactly \(k\) packed passengers. In particular,
\begin{equation}\label{eq:Hspecial}
 H_n(1)=2^n,\qquad
 H_n(0)=\#\{s\in\{0,1\}^n:p(s)=0\}.
\end{equation}
For example,
\begin{equation}\label{eq:H6}
 H_6(z)=44+12z+5z^2+2z^3+z^4,
\end{equation}
whose coefficients sum to \(64\).

\subsection{A rational generating function}

\begin{theorem}\label{thm:Hgf}
The unrestricted generating function is
\begin{equation}\label{eq:Hgf}
 \Hcal(x,z)=
 \frac{1+(1-z)x+(1-z)x^2}
 {1-(1+z)x-(1-z)x^2-(1-z)x^3}.
\end{equation}
\end{theorem}

\begin{proof}
Writing \(n=m+q\) and using \eqref{eq:Fq}, we obtain
\[
 \Hcal(x,z)=\sum_{q\geq0}x^q\Fcal_q(x,z)
 =\sum_{q\geq0}x^qB(x,z)^{q+1}
 =\frac{B(x,z)}{1-xB(x,z)}.
\]
The geometric-series expansion is valid formally because \(xB(x,z)\)
has zero constant term in \(x\). Substitution of \eqref{eq:B} and
simplification give \eqref{eq:Hgf}.
\end{proof}

\begin{corollary}\label{cor:Hrec}
For \(n\geq3\),
\begin{equation}\label{eq:Hrec}
 H_n(z)=(1+z)H_{n-1}(z)
       +(1-z)H_{n-2}(z)+(1-z)H_{n-3}(z),
\end{equation}
with \(H_0(z)=1\), \(H_1(z)=2\), and \(H_2(z)=4\).
\end{corollary}

\begin{proof}
Multiply \eqref{eq:Hgf} by its denominator and compare coefficients of
\(x^n\). The numerator has degree \(2\), so the homogeneous recurrence
starts at \(n=3\). The result is this interesting Tribonacci sequence.
\end{proof}

\subsection{A characteristic-root formula}

For a fixed \(z\in\mathbb{C}\), let
\begin{equation}\label{eq:charpoly}
 f_z(r)=r^3-(1+z)r^2-(1-z)r-(1-z).
\end{equation}

\begin{proposition}\label{prop:Hexplicit}
If the roots \(r_1,r_2,r_3\) of \(f_z\) are distinct, then
\begin{equation}\label{eq:Hexplicit}
 H_n(z)=\sum_{j=1}^3
 \frac{r_j^2+(1-z)r_j+(1-z)}
      {3r_j^2-2(1+z)r_j-(1-z)}\,r_j^n,
 \qquad n\geq0.
\end{equation}
\end{proposition}

\begin{proof}
If the roots are distinct, none is zero: a zero root would force
\(z=1\), when \(f_1(r)=r^2(r-2)\). The denominator of
\eqref{eq:Hgf} factors as \(\prod_j(1-r_jx)\). Writing its partial
fraction decomposition as \(\sum_j C_j/(1-r_jx)\), evaluation at
\(x=1/r_j\) gives
\[
 C_j=\frac{r_j^2+(1-z)r_j+(1-z)}{f_z'(r_j)}.
\]
Extracting the coefficient of \(x^n\) proves the formula.
\end{proof}

At values of \(z\) for which roots coincide, the combined expression in
\eqref{eq:Hexplicit} is interpreted by continuity, rather than by
substitution into its individual summands. Such a limit exists because
each \(H_n(z)\) is a polynomial. The recurrence remains valid at every
specialization; for instance, \(H_n(1)=2^n\).

\section{Probability and limiting behavior}
\label{sec:probability}

\subsection{Uniform and fixed-count models}

Let \(s\) be uniformly distributed on \(\{0,1\}^n\). Equivalently, its
coordinates are independent Bernoulli random variables with parameter
\(1/2\). Then
\begin{equation}\label{eq:probH}
 \Prob(p(s)=k)=\frac{[z^k]H_n(z)}{2^n},\qquad
 \E[z^{p(s)}]=\frac{H_n(z)}{2^n}.
\end{equation}
Conditioning on \(m(s)=m\) gives the uniform distribution on \(S_{m,n}\),
and hence
\begin{equation}\label{eq:probG}
 \Prob(p(s)=k\mid m(s)=m)=
 \frac{[z^k]G_{m,n}(z)}{\binom nm}.
\end{equation}

\begin{proposition}\label{prop:means}
For \(n\geq3\),
\begin{align}
 \E[p(s)]&=\frac{n-2}{8},\label{eq:mean}\\
 \E[p(s)\mid m(s)=m]&=
 \frac{m(m-1)(m-2)}{n(n-1)},\qquad 0\leq m\leq n.
 \label{eq:meancond}
\end{align}
For \(n<3\), both expectations are zero.
\end{proposition}

\begin{proof}
There are \(n-2\) possible centers of occupied triples. Under the
unrestricted model, each corresponding triple is occupied with
probability \(1/8\). Under the fixed-count model, that probability is
\[
 \frac{\binom{n-3}{m-3}}{\binom nm}
 =\frac{m(m-1)(m-2)}{n(n-1)(n-2)}
\]
when \(m\geq3\), and is zero otherwise. Linearity of expectation gives
both formulas. Independence of overlapping triples is not required.
\end{proof}

\subsection{The probability of no packed passengers}

Set \(T_n=H_n(0)\) and \(Q_n=T_n/2^n\). Specializing \eqref{eq:Hrec}
at \(z=0\) gives
\begin{equation}\label{eq:trib}
 T_n=T_{n-1}+T_{n-2}+T_{n-3}\quad(n\geq3),\qquad
 T_0=1,\quad T_1=2,\quad T_2=4.
\end{equation}
These numbers count binary strings avoiding \(111\).

\begin{theorem}\label{thm:Qlimit}
Let \(\lambda\) be the real root of \(r^3-r^2-r-1=0\). Then
\begin{equation}\label{eq:Qasymptotic}
 Q_n\sim C\left(\frac{\lambda}{2}\right)^n,\qquad
 C=\frac{\lambda^2+\lambda+1}{3\lambda^2-2\lambda-1},
\end{equation}
where \(\lambda\approx1.839286755\) and \(C\approx1.13745\).
In particular, \(\lim_{n\to\infty}Q_n\longrightarrow0\).
\end{theorem}

\begin{proof}
The cubic \(g(r)=r^3-r^2-r-1\) has discriminant \(-44\), so it has
one real root and two nonreal conjugate roots \(\alpha,\overline\alpha\).
Since \(g(1)<0<g(2)\), its real root satisfies \(1<\lambda<2\).
The product of the roots is \(1\), so
\(\lvert\alpha\rvert=\lambda^{-1/2}<1\).
Applying Proposition~\ref{prop:Hexplicit} at \(z=0\), we obtain
\[
 T_n=C\lambda^n+D\alpha^n+\overline D\,\overline\alpha^{\,n}
\]
for some \(D\in\mathbb C\), with \(C\) as in
\eqref{eq:Qasymptotic}. Both its numerator and denominator are positive:
the latter is \((3\lambda+1)(\lambda-1)\). Thus \(C>0\), and
\(T_n\sim C\lambda^n\). Dividing by \(2^n\) proves the claim.
\end{proof}

The limit alone also has an elementary proof. Partition the first
\(3\lfloor n/3\rfloor\) seats into disjoint triples. If no packed
passenger occurs, none of these disjoint triples is \(111\).
Their independence therefore gives
\[
 0\leq Q_n\leq(7/8)^{\lfloor n/3\rfloor}\longrightarrow0.
\]

\section{Full strings and Fibonacci numbers}
\label{sec:full}

\begin{definition}[Full row]\label{def:full}
A binary string is full if changing any one of its zero entries to \(1\)
makes the newly occupied seat packed. For \(n\geq1\), let \(a_n\) be
the number of full strings of length \(n\). An all-ones string is full
vacuously.
\end{definition}

This is a condition on empty seats, not a requirement that every
currently occupied seat be packed. Nor does it require that the original
string have packed count zero.

For \(n\geq2\), the zero positions of a full row form an independent
set of the path on the \(n-2\) interior positions, and conversely.
Consequently, the Fibonacci and fixed-cardinality formulas below
are the usual path-independent-set counts; compare
Codara and D'Antona \cite[Lemma~3 and Remark~4]{CD}.\footnote{Full text:
\url{https://arxiv.org/html/1211.2251}.}
We include direct seating proofs for completeness.

\begin{proposition}\label{prop:fullchar}
For \(n\geq2\), a string is full if and only if its first and last entries
are \(1\) and it contains no \(00\). For \(n=1\), only the string \(1\)
is full.
\end{proposition}

\begin{proof}
An endpoint zero cannot become packed. An interior zero becomes packed
precisely when both of its neighbors are \(1\). Thus every zero must
be interior and isolated, which is equivalent to the stated conditions.
\end{proof}

\begin{theorem}\label{thm:afib}
Let \(F_0=0\), \(F_1=1\), and \(F_n=F_{n-1}+F_{n-2}\) for \(n\geq2\).
Then \(a_n=F_n\) for \(n\geq1\).
\end{theorem}

\begin{proof}
The initial counts are \(a_1=a_2=1\). For \(n\geq3\), a full string
ends either in \(11\) or in \(01\). In the first case, deleting the last
\(1\) gives a full string of length \(n-1\). In the second case, the last
three entries are \(101\), and deleting the final \(01\) gives a full
string of length \(n-2\). These operations are bijections, reversed by
appending \(1\) and \(01\), respectively. Hence
\(a_n=a_{n-1}+a_{n-2}\), proving the result.
\end{proof}

Under the unrestricted uniform model,
\begin{equation}\label{eq:fullprob}
 \Prob(s\text{ is full})=\frac{F_n}{2^n}
 \sim\frac{1}{\sqrt5}\left(\frac{\varphi}{2}\right)^n
 \longrightarrow0,\qquad \varphi=\frac{1+\sqrt5}{2}.
\end{equation}
Indeed, with \(\psi=(1-\sqrt5)/2\), the exact ratio is
\(\bigl((\varphi/2)^n-(\psi/2)^n\bigr)/\sqrt5\).
The definition of \(a_n\) is restricted to \(n\geq1\); the vacuously
full empty string is not included in the identity \(a_n=F_n\).

\begin{proposition}\label{prop:fullfixed}
For \(1\leq m\leq n\), the number of full strings in \(S_{m,n}\) is
\begin{equation}\label{eq:fullfixed}
 \binom{m-1}{n-m}.
\end{equation}
Here \(\binom ab=0\) when \(a\geq0\) and either \(b<0\) or \(b>a\).
There are no full strings with \(m=0\) and \(n\geq1\).
\end{proposition}

\begin{proof}
Start with \(m\) ones. Every zero must occupy a distinct one of the
\(m-1\) internal gaps between consecutive ones, and every such insertion
produces a full string. Choose the \(n-m\) gaps.
\end{proof}

It follows that, for \(1\leq m\leq n\),
\begin{equation}\label{eq:fullfixedprob}
 \Prob(s\text{ is full}\mid m(s)=m)
 =\frac{\binom{m-1}{n-m}}{\binom nm}.
\end{equation}
Summing \eqref{eq:fullfixed} over \(m\) gives
\begin{equation}\label{eq:fibidentity}
 F_n=\sum_{m=1}^n\binom{m-1}{n-m},\qquad n\geq1.
\end{equation}

\section{Rectangular seating arrangements}
\label{sec:grid}

Let \(n,r\geq1\), and let \(s=(s_{i,j})\) be a binary \(n\times r\)
matrix. A seat is packed when it and its four orthogonally adjacent
seats are occupied. Boundary seats are never packed: having all
neighbors that happen to exist occupied is not sufficient.
The packed count is
\begin{equation}\label{eq:p2}
 p_2(s)=\sum_{i=2}^{n-1}\sum_{j=2}^{r-1}
 s_{i,j}s_{i-1,j}s_{i+1,j}s_{i,j-1}s_{i,j+1}.
\end{equation}
In particular, \(p_2(s)=0\) whenever \(n<3\) or \(r<3\).
Define
\begin{equation}\label{eq:H2def}
 H^{(2)}_{n,r}(z)=
 \sum_{s\in\{0,1\}^{n\times r}}z^{p_2(s)}.
\end{equation}
For square grids, write \(H^{(2)}_n(z)=H^{(2)}_{n,n}(z)\).

Let \(I_{n,r}=\{2,\ldots,n-1\}\times\{2,\ldots,r-1\}\), with either
factor understood as empty when its upper bound is below \(2\).
For \(v=(i,j)\in I_{n,r}\), set
\[
 N(v)=\{(i,j),(i-1,j),(i+1,j),(i,j-1),(i,j+1)\}.
\]
For \(A\subseteq I_{n,r}\), write
\[
 V(A)=\bigcup_{v\in A}N(v),\qquad V(\varnothing)=\varnothing.
\]
We use \(V\) here to distinguish neighborhood unions from the isolated-seat
statistic \(U(s)\) in Section~\ref{sec:minimum}.

\begin{theorem}[Neighborhood-union expansion]\label{thm:grid}
For \(n,r\geq1\),
\begin{equation}\label{eq:H2explicit}
 H^{(2)}_{n,r}(z)=
 \sum_{A\subseteq I_{n,r}}2^{nr-|V(A)|}(z-1)^{|A|}.
\end{equation}
Consequently, for \(0\leq k\leq |I_{n,r}|\),
\begin{equation}\label{eq:H2coef}
 [z^k]H^{(2)}_{n,r}(z)=
 \sum_{\substack{A\subseteq I_{n,r}\\|A|\geq k}}
 (-1)^{|A|-k}\binom{|A|}{k}\,2^{nr-|V(A)|}.
\end{equation}
\end{theorem}

\begin{proof}
For each \(v\in I_{n,r}\), let \(X_v(s)=\prod_{u\in N(v)}s_u\).
Since \(X_v(s)\in\{0,1\}\),
\[
 z^{p_2(s)}
 =\prod_{v\in I_{n,r}}\bigl(1+(z-1)X_v(s)\bigr)
 =\sum_{A\subseteq I_{n,r}}(z-1)^{|A|}
       \prod_{v\in A}X_v(s).
\]
For a fixed \(A\), the final product equals \(1\) precisely when every
entry indexed by \(V(A)\) is \(1\). The other \(nr-|V(A)|\) entries are
free, giving \(2^{nr-|V(A)|}\) matrices. Summing over \(s\) proves
\eqref{eq:H2explicit}. Expanding each \((z-1)^{|A|}\) by the binomial
theorem and extracting \([z^k]\) gives \eqref{eq:H2coef}.
\end{proof}

The expansion also covers grids with no interior positions:
only \(A=\varnothing\) contributes, giving \(H^{(2)}_{n,r}(z)=2^{nr}\).
Although exact, the formula requires \(2^{|I_{n,r}|}\) subsets in its
direct form; it is not a polynomial-time coefficient algorithm.

Theorem~\ref{thm:grid} is also a specialization of a general
hypergraph-polynomial identity. For a hypergraph \(\mathcal K\) on
\(N\) vertices, put
\[
 P_{\mathcal K}(u,z)=
 \sum_{W\subseteq V(\mathcal K)}
 u^{|W|}z^{|\{e\in E(\mathcal K):e\subseteq W\}|},
 \qquad
 S_{\mathcal K}(a,b)=
 \sum_{A\subseteq E(\mathcal K)}
 a^{|\bigcup_{e\in A}e|}b^{|A|}.
\]
Tadesse's transformation \cite[Corollary~2.2]{T}
is\footnote{Full text:
\url{https://arxiv.org/pdf/1303.4520}.}
\[
 P_{\mathcal K}(u,z)=(1+u)^N
 S_{\mathcal K}\left(\frac{u}{1+u},z-1\right).
\]
Take the grid cells as vertices, the neighborhoods \(N(v)\) as
hyperedges, and \(u=1\). Then \(N=nr\), the left side is
\(H^{(2)}_{n,r}(z)\), and the right side is exactly
\eqref{eq:H2explicit}. Thus the grid expansion is included as an
application of this transformation, which provides an yet another interesting connection to hypergraphs.

\begin{example}[A \(4\times4\) grid]\label{ex:grid4}
There are four interior positions. Their neighborhood unions give the
following sums:
\[
 \begin{array}{cr}
 \toprule
 j & \displaystyle\sum_{\substack{A\subseteq I_{4,4}\\|A|=j}}
                    2^{16-|V(A)|}\\
 \midrule
 0 & 65536\\
 1 & 8192\\
 2 & 1536\\
 3 & 256\\
 4 & 16\\
 \bottomrule
 \end{array}
\]
For a singleton, the union has size \(5\). Every pair of distinct
interior neighborhoods has union size \(8\); every triple has union
size \(10\); and all four together have union size \(12\).
Thus the five entries are, respectively,
\(2^{16}\), \(4\cdot2^{11}\), \(6\cdot2^8\), \(4\cdot2^6\), and
\(2^4\).
Substitution into \eqref{eq:H2explicit} gives
\begin{align*}
 H^{(2)}_4(z)
 &=65536+8192(z-1)+1536(z-1)^2\\
 &\hspace{1.5em}+256(z-1)^3+16(z-1)^4\\
 &=58640+5824z+864z^2+192z^3+16z^4.
\end{align*}
In particular, exactly \(864\) of the \(4\times4\) matrices have two
packed passengers, and \(H^{(2)}_4(1)=65536=2^{16}\).
\end{example}

\end{document}